\documentclass[12pt]{article}

\title{Computing the rational torsion of an elliptic curve using
Tate normal form}
\author{
Irene Garc\'{\i}a $\;$ Miguel A. Olalla\thanks{Partially supported by
Junta de Andaluc\'{\i}a, Ayuda a Grupos FQM 218} $\;$ Jos\'e M.
Tornero$^\ast$
\\ Facultad de Matem\'aticas, Universidad de Sevilla \\ Apdo. 1160
E-41080 SEVILLA (SPAIN)}
\date{Email: olalla@algebra.us.es, tornero@algebra.us.es}

\begin{document}

\maketitle

{\bf Abstract.} It is a classical result (apparently due to Tate)
that all elliptic curves with a torsion point of order $n$ ($4
\leq n \leq 10$, or $n = 12$) lie in a one-parameter family.
However, this fact does not appear to have been used ever for
computing the torsion of an elliptic curve. We present here a
extremely down--to--earth algorithm using the existence of such a
family.


\section*{1. Tate and Weierstrass normal forms}

An elliptic curve is a plane smooth affine (respectively
projective) curve defined by a cubic (homogeneous) polynomial. All
these curves are known to be birrationally equivalent (that is,
isomorphic as algebraic varieties, up to a finite number of
points) to one which equation has the form $Y^2 + a_1 XY + a_3Y =
X^3 + a_2 X^2 + a_4 X + a_6$. When all the coefficients lie in a
field $K$, the set of points in the curve with both coordinates in
$K$ admits a group structure (\cite{Cassels1}, \cite{Cassels},
\cite{Husemoller}) with the inner operation defined by the
classical chord--tangent procedure. This group is then noted
$E(K)$. For historical reasons we will note this operation
additively and so we will write $2P$ for $P+P$. As the unit
element is usually taken to be the only point at infinity (say
${\cal O}$), we can restrict ourselves to affine points.

The Mordell--Weil theorem states that, if $K$ is a number field,
$E(K)$ is always a finitely generated abelian group
(\cite{Cassels}, \cite{Husemoller}). The torsion subset of $E(K)$
is hence a finite subgroup, noted $E_T(K)$. The strongest result
concerning $E_T ({\mathbf Q})$ is due to B. Mazur and explicitly
states all groups which can appear as torsion subgroups of
elliptic curves defined over ${\mathbf Q}$:

\vspace{.3cm}

{\bf Theorem (Mazur).--} (\cite{Mazur}, \cite{Mazur2}) Let $E$ be
an elliptic curve defined over ${\mathbf Q}$. Then its torsion
group $E_T ({\mathbf Q})$ is either isomorphic to $C_n$ (the
cyclic subgroup of $n$ elements) for $n=1,2,..., 10,12$ or to $C_2
\times C_{2n}$ for $n=1,2,3,4$. All of these possibilities
actually occur.

\vspace{.3cm}

The aim of this paper is giving an efficient procedure, different
from the usual ones, still very lowbrow, for computing the torsion
subgroup of an elliptic curve defined over the rationals. First of
all we must put the curve into a more manageable form.

For a general elliptic curve it is known (see, for instance,
\cite{Cassels1}, \cite{Cassels}, \cite{Husemoller}) that using
linear changes of variables, one can take the equation defining the
elliptic curve into an easier one of the type $Y^2 = X^3 + AX +
B$. This is known as Weierstrass (short) normal form.

A straightforward computation proves that the only linear changes
of variables preserving Weierstrass normal form are those given by
$$
\left\{ \begin{array}{ccc}
X & \longmapsto & u^2 X' \\
Y & \longmapsto & u^3 Y'
\end{array} \right.
$$
for some $u \in {\bf Q}$. Such a change takes the curve defined by
$Y^2 = X^3 + AX + B$ into the one defined by $Y^2 = X^3 + (A/u^4)X
+ (B/u^6)$. This argument shows that one can always assume $A$
and $B$ to be in ${\mathbf Z}$. It also implies that the number
$A^3/B^2$ is an invariant of the equivalence class of elliptic
curves in Weierstrass form up to linear changes of variables.

Of course, even if two curves $Y^2 = X^3 + AX + B$ and
$Y^2 = X^3 + CX + D$ verify $A^3/B^2 = C^3/D^2$ this does not mean
they are equal up to some linear change of variables of the
previous form. In fact, it is fairly elementary proving that this
happens if and only if the following condition hold: there exists
a rational solution $u$ for the system
$$
\left\{ u^4 = \frac{A}{C}, \; \; u^6 = \frac{B}{D} \right\},
$$
with the obvious arrangements for the cases in which any of the
coefficients vanishes.

\vspace{.3cm}

In addition, if the curve is already known to have one rational
point of order $n>3$, one can choose to put the equation of the
curve in the form $Y^2 + bXY + cY = X^3 + dX^2$, also using
nothing but linear changes of variables. This second formula is
called Tate normal form (\cite{Husemoller}).

\section*{2. The Lutz -- Nagell theorem}

Most classical algorithms for computing rational torsion of
elliptic curves are based on the following result, achieved
independently by Lutz and Nagell (\cite{Nagell},
\cite{Lutz}):

\vspace{.3cm}

{\bf Theorem (Lutz -- Nagell).--} Let $E$ be an elliptic curve
defined over ${\mathbf Q}$, given by a Weierstrass equation $Y^2 = X^3
+ AX + B$ with $A, B \in {\bf Z}$, and let $P = (\alpha, \beta)
\in E_T ({\mathbf Q})$. Then

\begin{enumerate}
\item[(a)] Both $\alpha$ and $\beta$ are in ${\bf Z}$.

\item[(b)] Either $\beta = 0$ or $\beta^2 | (4A^3 + 27B^2)$.
\end{enumerate}

\vspace{.3cm}

Clearly $\beta = 0$ is equivalent to $2P = O$ and this 2-torsion
part is rapidly computable. For computing the remaining points (if
there are any) we simply factorize $\Delta = 4A^3 + 27 B^2$. This
quantity, called the discriminant of $E$, will be most important
in the sequel. For every square divisor, say $m^2$ of $\Delta$, we
compute the integral solutions to $X^3 + AX + (B-m^2)$. If we
actually find an integral root, say $n$, we only have to check
whether $(n,m)$ is a torsion point, which only involves
computing, at most, $12(n,m)$, following Mazur's Theorem.

Simple as it is, this algorithm is not very efficient, being its
major drawback the necessity of factoring $\Delta$. This is the
algorithm presented, for instance, in \cite{Cohen} and
\cite{Cremona}.

\section*{3. Good reduction: a first bound}

The first step in our algorithm will be a reasonable bound for the
size of $E_T ({\mathbf Q})$. The existence of the group structure
in an elliptic curve does not depend on the field we are taking
coordinates in. So, for instance, if $A,B \in {\mathbf Z}$, then
for all primes $p$, the same equation defining an elliptic curve
over ${\mathbf Q}$ defines an elliptic curve over the finite field
${\mathbf F}_p$. The relationship between these two curves can
help us in our purpose, using the next result (\cite{Cassels},
\cite{Husemoller}):

\vspace{.3cm}

{\bf Theorem.--} Let $E$ be an elliptic curve in Weierstrass form
$Y^2 = X^3 + AX + B$, with $A,B \in {\mathbf Z}$. If $p>2$ is a
prime number such that it does not divide $\Delta$, then the
mapping
\begin{eqnarray*}
\mbox{red}_p: E_T ({\mathbf Q})& \longrightarrow & E ({\mathbf F}_p) \\
(\alpha_1, \alpha_2) & \longmapsto & (\overline{\alpha_1},
\overline{\alpha_2}) \\
{\mathcal O} & \longmapsto & {\mathcal O}
\end{eqnarray*}
is an injective group homomorphism (where $\overline{\alpha_i}$ denotes
the residue classes of $\alpha_i$ modulo $p$).

\vspace{.3cm}

Primes which do not divide $\Delta$ are called good primes and the
induced group homomorphisms are called good reductions. So,
choosing some prime $p$ not dividing $\Delta$ and computing how
many points lie in $E({\mathbf F}_p)$ we must obtain a multiple of
the order of $E_T({\mathbf Q})$. Our practical choice has been
taking three primes (as small as possible), computing the number
of points in each case and finding the greatest common divisor of
all those quantities. In most cases, this bound was found to be
the actual order of $E_T({\mathbf Q})$.

\vspace{.3cm}

There are, however, some cases which does not fit this scheme. For
example, the curve defined by $Y^2 = X^3 + X$ has the property
that, $E_T({\mathbf Q}) = C_2$ but, for every good prime $p$ the
order of $E({\mathbf F}_p)$ is divisible by 4.

Trying then to be a bit more accurate, we computed not only the
order of $E({\mathbf F}_p)$, but also how many elements of order 2
it had. So, if $E({\mathbf F}_p)$ presented more points of order 2
than $E$ itself, our choosing for the bound can be smaller than
the order of $E({\mathbf F}_p)$. In the above example, as $E_3$ is
isomorphic to $C_4$ and $E ({\mathbf F}_5)$ is isomorphic to $C_2
\times C_2$, the bound actually found is the order of the group
$E_T({\mathbf Q})$.

So, if $|E({\mathbf F}_p)| = M$, the number of points of order 2
in $E$ is $s$ and the number of points of order 2 in $E({\mathbf
F}_p)$ is $t$, the choosing of the bound goes like this:

\begin{enumerate}
\item[(a)] If $s=t$, we choose $M$.

\item[(b)] If $(s,t) \in \{ (0,1),(1,3) \}$ then we choose $M/2$.

\item[(c)] If $(s,t) = (0,3)$ then we can choose $M/4$ as the bound.
\end{enumerate}

Note that one needs the fact that $E({\mathbf F}_p)$ is a finite
group with, at most, three elements of order 2.

\section*{4. Points of given order}

We will explain now how to decide when an elliptic curve defined
over the rationals has a point of a given order, say $n$, where
$n=4,...,10,12$. First we need a result on parametrization of
torsion structures. Most cases are proved (quite
straightforwardly) in \cite{Husemoller}. Also see \cite{Kubert}
for a more exhaustive table, without any proofs.

\vspace{.3cm}

{\bf Theorem.--} Every elliptic curve with a
point $P$ of order $n=4,...,9,10,12$ can be written in the
following Tate normal form
$$
Y^2 + (1-c) XY  - b Y = X^3 - bX^2,
$$
with the following relations:
\begin{itemize}
\item[(1)] If $n=4$, $b = \alpha, \; c =0$.
\item[(2)] If $n=5$, $b = \alpha, \; c = \alpha$.
\item[(3)] If $n=6$, $b = \alpha + \alpha^2, \; c = \alpha$.
\item[(4)] If $n=7$, $b = \alpha^3 - \alpha^2, \; c = \alpha^2
- \alpha$.
\item[(5)] If $n=8$, $b = (2\alpha - 1)(\alpha - 1), \;
c= b/\alpha$.
\item[(6)] If $n=9$, $c=\alpha^2 (\alpha-1), \; b = c(\alpha
(\alpha-1)+1)$.
\item[(7)] If $n=10$, $c= (2\alpha^3-3\alpha^2+\alpha)/\left[
\alpha - (\alpha -1)^2 \right], \; b=c\alpha^2/\left[ \alpha-
(\alpha-1)^2 \right]$.
\item[(8)] If $n=12$, $c= (3\alpha^2-3\alpha+1)(\alpha-
2\alpha^2)/ (\alpha-1)^3, \; b= c(2\alpha-2\alpha^2-1)/(\alpha-1)$.
\end{itemize}

\vspace{.3cm}

Suppose then that we want to check if a given curve $E$ defined by
$Y^2 = X^3 + AX + B$ has a point of order $n$. Assume it posseses such a
point: therefore $E$ must be isomorphic to one curve lying in the
one--parameter family. Then we simply compute the Weierstrass
normal form of a generic curve in the familiy and check the
conditions given at the end of section 1 for two curves in
Weirestrass form to be isomorphic.

\vspace{.3cm}

{\bf Example.--} Let us give an example with $n=5$. Suppose that
we would like to know if our curve $Y^2=X^3+12933X-2285226$ (this
is curve 110A1(C) from \cite{Cremona})  has a point of order 5. If
it is the case, the curve must be isomorphic, by a linear change
of variables, to one lying in the family
$$
Y^2 + (1-\alpha)XY - \alpha Y = X^3 - \alpha X^2.
$$

So, taking this general equation to Weierstrass form we obtain
an equation which we will note $Y^2 = X^3 + A_5(\alpha)X + B_5(\alpha)$.
Should this curve be isomorphic to ours, it must hold
$$
\frac{A_5(\alpha)^3}{B_5(\alpha)^2} = \frac{12933^3}{2285226^2},
$$
which sums up to an equation in the variable $\alpha$ (in our case, of
degree 12).

This equation will be called the final polynomial for $n=5$. For
every root $\alpha_0$ we have to check if there is some $u \in {\bf Q}$
verifying
$$
\left\{ u^4 = \frac{A}{A_5(\alpha_0)}, \; \; u^6 = \frac{B}{B_5(\alpha_0)}
\right\}.
$$

If there is then we have a point of order 5, which is easily
calculated, as $(0,0)$ is a point of order 5 in the Tate normal form. If not,
then there are no points of order 5 in
 $E$.

In our example, the only roots were $-1/10$ and 10. Besides,
$$
A_5(10) = A, \; \; B_5 (10) = B,
$$
so in fact there is a point of order 5 in our curve. Tracing back
the changes of variables a point of order 5 turns out to be $(123,1080)$.

\vspace{.3cm}

The only remaining case is $n=3$ that is, we need a procedure for
deciding if an elliptic curve has a point of order 3. There is
also a Tate normal form for this case, but it has some
inconveniences, being the heaviest one that the family of curves
depends now on two parameters. However, there is a well-known
property which can be used (\cite{Cassels1}):

\vspace{.3cm}

{\bf Proposition.--} Let $E$ be an elliptic curve given by a
Weierstrass equation $Y^2 = X^3 + AX + B$. Then $E$ has a point
$P$ of order 3 if and only if there is an integral solution to the
equation
$$
3X^4 + 6AX^2 + 12BX - A^2 =0.
$$

In this case, the solution is the first coordinate of $P$. In fact,
in the cited article one can find polynomials which characterize
points of any order. These polynomials become more complicated as
the order grows, but they also allow to obtain a obvious procedure
for deciding if there is any point of given order.

\section*{5. The algorithm}

Given an elliptic curve in Weierstrass form $Y^2 = X^3 +AX + B$,
in order to find its torsion group we proceed as follows:

\vspace{.3cm}

{\bf Step 1.} Compute the number of points with order 2, that is,
the rational solutions for $X^3 + AX + B$.

{\bf Step 2.} Pick the smaller five (for instance) good primes for
$E$ and compute a bound $M$ for the torsion as explained above.

{\bf Step 3.} If the number of rational solutions is either 0 or
1, then for every divisor $d$ of $M$, apply the procedure
described in the previous section to check if there is a point of
order $d$. If this is done is decreasing order, the first
affirmative answer gives us the group (which should be isomorphic
to $C_d$) and one generator: either the point which comes from
point $(0,0)$ in Tate normal form for $n=4,...,10,12$ or the point
directly obtained for $n=3$.

{\bf Step 4.} If the number of rational solutions is 3, then apply
the same procedure as above for every divisor $d$ of $M/2$. Now
the first affirmative answer gives us the group (which must be
$C_2 \times C_d$) and a set of generators (the points of order 2
and the point which comes from point $(0,0)$ in Tate normal form).

\section*{6. Explicit calculations}

In this section, we will show the computations that led us to the
implementation of our algorithm in Maple, currently available by
anonymous ftp at
 \verb|ftp://alg7.us.es/pub/Programs/| (comments
in Spanish so far...).

So we fix an elliptic curve $E$, given by $Y^2 = X^3 + AX + B$
with $A,B \in {\mathbf Z}$ and we want to know if there is a point of
order $n$ on it. For all cases (except $n=3$) we know this implies
solving an equation on a parameter $\alpha$ which comes from
the parametrizations of Tate normal form.

However, one may find that ``classical'' parametrizations, though
the simplest ones, are not necessarily the most convenient for our
purpose. As we will need to compute the rational solutions of a
polynomial in ${\mathbf Z}[X]$, which the best parameter is depends
heavily on which root finding method is to be used.

Our choice was the algorithm developed in \cite{SIAM}, so we had to
take into account that the complexity of finding the rational roots a
polynomial in ${\mathbf Z}[X]$, say $f(X) = \sum a_iX^i$, of degree
$n$, is $O(\log^2 ||f||)$, where
$$
||f|| = \sum |a_i|,
$$
so one may choose a parameter which minimizes $||f||$ when $f$ is the
final polynomial. Such a parameter will be called a {\em minimal
parameter}.

\vspace{.3cm}

{\bf Case $n=4$.} We will do this in detail. The general equation
was
$$
Y^2 + XY - \alpha Y = X^3 - \alpha X^2,
$$
provided $b \neq 0, -1/16$.

Once it is taken to Weierstrass normal form, it sums up to
$$
Y^2 = X^3 + A_4( \alpha )X + B_4(\alpha),
$$
where
$$
A_4(\alpha) = -432\alpha^2 - 432\alpha - 27, \; \; B_4(\alpha) =
-3456\alpha^3 + 6480\alpha^2 + 1296\alpha + 54.
$$

So the final polynomial for $\alpha$,
$B(\alpha)^2A^3-A(\alpha)^3B^2$ results
$$
\begin{array}{rcl}
P_4 (\alpha) & = & 2^{12} 3^6  \Delta \alpha^6 -
2^{12} 3^7 \left( 5A^3-27B^2 \right) \alpha^5 \\
& & \quad + 2^8 3^7 \left( 59A^3+459B^2 \right) \alpha^4
+ 2^9 3^6  11  \Delta \alpha^3 \\
& & \quad + 2^4 3^7 17  \Delta \alpha^2 +
2^4 3^7 \Delta \alpha + 3^6 \Delta
\end{array}
$$

Our next step is then to find a minimal parameter (that is, a
parameter minimizing the norm of its final polynomial). So we find
a new parameter $\beta = r \alpha + s$. Obviously we need our new
final polynomial, $F_4 (\beta)$ to lie in ${\mathbf Z} [X]$ so it
is plain that the natural choosing for $r$ must be $1/12$. Then
we look for a rational $s$ which minimizes $||F_4||$. As $F_4$ was
to lie in ${\mathbf Z}[X]$ the possible denominators were bounded
(actually they had to be a divisor of 12). We find a minimum for
$s=1/12$ so we took $\alpha = (\beta+1)/12$ and
$$
\begin{array}{rcl}
F_4 (\beta) & = & \Delta \beta^6 - 6 \left(34A^3-135B^2 \right)
\beta^5 + 3\left( 851A^3+2646B^2 \right) \beta^4 \\
& & \; \; \; + 4 \left(313A^3+5940B^2 \right) \beta^3 -6 \left(
95A^3+2646B^2 \right)  \beta^2 \\
& & \; \; \; -24 \left( A^3 - 135B^2 \right) \beta + 49A^3 - 216 B^2.
\end{array}
$$

If we set $N = \max \{ |A|^3, |B|^2 \}$ then
$$
||F_4|| \leq 56667 N \simeq 2^8 3^2 5^2 N.
$$

\vspace{.3cm}

We present below all the minimal parameters along with bounds for
the seminorm of the final polynomials, calculated as above.

\vspace{.3cm}

{\bf Case $n=5$.} $\beta = \alpha$, $\mbox{deg} (F_5) = 12$,
$||F_5|| \leq 898312 N \simeq 2^{12} 3^2 5^2 N$.


\vspace{.3cm}

{\bf Case $n=6$.} $\alpha = \beta/3 - 1/3$, $\mbox{deg} (F_6) =
12$, $||F_6|| \leq 2220071 N \simeq 2^43^25^6 N$.

\vspace{.3cm}

{\bf Case $n=7$.} $\beta = \alpha$, $\mbox{deg}(F_7) = 18$,
$||F_7|| \leq 110725743 N \simeq  2^{22}3^{3}N$.

\vspace{.3cm}

{\bf Case $n=8$.} $\alpha = \beta+1$, $\mbox{deg} (F_8) = 24$,
$||F_8|| \leq  46702469380 N \simeq 2^{9} 3^{5} 5^8 N$.

\vspace{.3cm}

{\bf Case $n=9$.} $\beta = \alpha$, $\mbox{deg} (F_9) = 36$,
$||F_9|| \leq 11353024920 N \simeq 2^{10} 3^6 5^6 N$.

\vspace{.3cm}

Cases $n=10$ and $n=12$ can of course be worked out in the same
way but the polynomials get quite unpractical. As
$$
E_T({\mathbf Q}) = \left\{
\begin{array}{ccc}
C_{10} & \Longleftrightarrow & C_2,C_5 \subset E_T({\mathbf Q}) \\
\\
C_{12} & \Longleftrightarrow & C_4,C_6 \subset E_T({\mathbf Q})
\end{array} \right.
$$
there is no necessity of finding the actual polynomials $F_{10}$ and
$F_{12}$.
In these cases, the generator can be easily computed using
the duplication formula.

The leading coefficient of all final polynomials turns out to
be $\Delta$. Indeed, one can look for a parameter such that the leading
coefficient and the independent term of its final polynomials are $\Delta$.
So, if the factorization of $\Delta$ is known, this final polynomials can
speed up the process, as all the possible rational roots of the final
polynomials are known in advance.

\section*{7. Complexity and some examples}

As in the previous section, let
$$
N = \max \left\{ |A|^3, |B|^2 \right\}.
$$

We will show that the running time of our algorithm is ${\cal O}
(K\log^2 N)$ for some $K \in {\mathbf N}$. Unless otherwise
stated, \cite{Cohen} is the reference here for the details.

\vspace{.3cm}

The computation of the points of order two can be clearly
accomplished in the expected time, using, for instance, the algorithm
given in \cite{SIAM}. Note that, should this be the case, it can
also be used for checking the existence of points with order three,
with the desired complexity.

\vspace{.3cm}

The bounding of the torsion consists only on arithmetical
operations on affine planes ${\mathbf F}_p$, with $p$ not dividing
$\Delta$. It is clear that there are primes smaller than $N$ which
not divide $\Delta$. Of course, it is known that arithmetical
operations with data bounded by $\Delta$ can be carried out in ${\cal
O} (\log^2 N)$ time.

\vspace{.3cm}

So it only remains checking step 3 (step 4 is analogous) for the
cases $n=4,...,10,12$. But note that all the coefficients of our
minimal polynomials are bounded by $cN$, for some natural $c$. This means
that, for a
 rational root, written in irreducible form $\alpha_0 =
\beta_0 /
 \gamma_0$, we have
$$
\left| \beta_0 \right|, \left| \gamma_0 \right| <  cN.
$$

Therefore, if we want to find out if there exists some $u \in {\bf
Q}$ such that $u^4 = A/A_n(\alpha_0)$ we only have to put $A /A_n
(\alpha_0)$ in irreducible form (that amounts to find the gcd and
divide) and compute the square root of its numerator and
denominator twice. All these operations can be carried out in the
expected time. If such an $u$ exists, it is just a matter of
arithmetical checking seeing if $u^6 = B/B_n(\alpha_0)$.

\vspace{.3cm}

Some time results of our algorithm are given in the following
examples table, using our MapleV routine.
$$
\begin{array}{rl}
E_1: & Y^2 = X^3 - 98D6E49C45C901B \cdot X + B5D1E097F653622F55B036\\
E_2: & Y^2 = X^3 - (A_2/A'_2)X + (B_2/B'_2) \\

E_3: & Y^2 = X^3 - (A_3/16)X - (B_3/32) \\
\end{array}
$$
where
$$
\begin{array}{cl}
A_2 = & \mbox{83ACFBAEC1BB1AC8EA33B897FDE9672AB898D04622635}/ \\
& \mbox{198248803F6F6429EC185BB2AB6D5DAE2C41BA0EC07AD5}/ \\
& \mbox{46CFF23FA458FCB36D8E85877CF0} \\

A'_2 = & \mbox{4E07B196F78B523E2BF8B93D9FF09BFF22E07284643617}/ \\
& \mbox{AD603BDEBE49E967484527B634E2990C1E19261C903}/ \\ &
\mbox{AAC97D0F23EE86534D5011DF9A71} \\
\end{array}
$$

$$
\begin{array}{cl}
B_2 = & \mbox{1594F960645253D0B7F933BFD50446DC3FC067CAFEDB11}/ \\
& \mbox{E76E7EBBDA0FCB2EF4AC34672D4B6469AD156134B7DEA}/ \\ &
\mbox{2FC9C7EFA07084E7695B18DBE22D436EEE2BB5EA14C26}/ \\ &
\mbox{D67AB385078CB862970A2B56D62C837D4E00A097490} \\ B'_2 = &
\mbox{AC51A232098DD799F2D035E3B630C2EE79B9C00C70B9013}/ \\ &
\mbox{071A6A0011C7A689A577D55A9BCCDA3FDCCE2FB25958A}/ \\ &
\mbox{D9D1F26D9D0D118651B0B5548FF001466E8D0BF7946D23F}/ \\ &
\mbox{9319CE52A96C7C9B2D0E37DEC87027D90109} \\ A_3 = &
\mbox{AF06EC915A7BC47C45CFBFA797633ACE67A79F7B381D29}/ \\ &
\mbox{BCCA243AABA230AF5BAD1058D41582134BECEF3F8DBB} \\ B_3 = &
\mbox{1BDA1A8FE9A5108EA7DB7FB6AE8EB3F7AB45A8D22614B}/ \\ &
\mbox{93FDB39D03E0B8324128145C706768EF5EE5BE37E68F4C} \\ &
\mbox{B5BC9EE31CC5B7EDA2C668D5CF0EFE9AA31F0B460EEB}
\end{array}
$$

\vspace{.3cm}

$$
\begin{array}{cccc}
\hline \hline
\mbox{Curve} & \mbox{Group} & \mbox{Generator} & \mbox{Time} \\
\hline \hline
E_1 & \; C_4 \; & \; (1C8CFC03,100F4DC00) \; & \; 2.33 s \\
E_2 & C_5 & (\lambda_1/\lambda_3^2,\lambda_2/\lambda_3^3) & \; 10.66s \\
E_3 & C_2 \times C_4 & (\lambda_4,\lambda_5) & \; 3.85s \\
\hline
\end{array}
$$
where
$$
\begin{array}{rl}
\lambda_1 = & \mbox{-1A8019538D071D5BFD9EEBA7B19BE9124EB6E592F0D15}/ \\
& \mbox{B0DD77D8016A58C} \\
\lambda_2 = & \mbox{-1626E05A34E5EA7E90A84BF3C4D604949BAA0DA532CDE1}/ \\
& \mbox{147804F9E6491E9E49F16F356882A85DA4C9785AC75C} \\
\lambda_3 = & \mbox{17C6E3C032F89045AD746684045E05} \\
\lambda_4 = & \mbox{-7A36225A2ADAAFA9B059FF46EE903619BD0C4E2AD3AA1}/4 \\
\lambda_5 = & \mbox{-897CE6A57036059EE6653F2FCC623CDCF4ADD7F02E202A}/8
\end{array}
$$

\vspace{.3cm}

The computations have been performed in a KMD300 computer. Note
 that our
current implementation does not include so far the root finding
 algorithm of
\cite{SIAM} but Maple V 5.1 built--in
routine, so it is hoped that a
complete implementation of our
algorithm will obtain even better results.

\vspace{.3cm}

We have compared our algorithm with, probably, the two most
efficient current ones: Pari/GP built--in procedure,
\verb|elltors| (see \cite{Pari}) and the routine \verb|Tor| from
the Maple package APECS (see \cite{Apecs}).

Pari/GP \verb|elltors| follows the algorithm described in
\cite{Doud}, using the analytic parametrization of the curve. It
is extremely fast and, besides, the periods of the lattice associated
to the curve are
directly computed by Pari/GP when you enter the curve
with the routine \verb|ellinit|. However, in some cases (we can not figure
out when or why),
 \verb|elltors| needs such a precision that it may become
unpractical.
 It remains, however, as our favourite choosing for medium--size
coefficients.
Here are the time results, expressed as (time for
\verb|ellinit|) + (time for \verb|elltors|), for the previous
examples, together with the precision (by 100) required.

$$
\begin{array}{ccc}
\hline \hline
\mbox{Curve} & \mbox{Precision} & \mbox{Time} \\
\hline \hline
E_1 & \; > 3600 \; & \; ?? \\
E_2 & \; 1300 & \; 1.05s \ + \ 11.92s \\
E_3 & \; 200 & \; 0.06s \ + \ 0.08s \\
\hline
\end{array}
$$

For $E_3$ \verb|elltors| gave an incorrect result: it output
$C_4$ for the structure. Hence there appears to be some minor
bug in the implementation. In all our computations, no errors
were found in \verb|elltors|  when working with cyclic groups.

APECS \verb|Tor| uses the polynomials mentioned at the end of
section 4. When you introduce a curve, which you must do before
computing its torsion, it computes a great deal of data, in
particular a bound for the torsion subgroup and other relevant
quantities. If data are moderately large (even significantly
smaller than the examples) this takes a huge lot of time: we mean
{\em hours} for the examples above. Anyway, its library is really
huge, so, for small--size coefficients, APECS will surely have
a lot of information (of course everything concerning rational torsion
points) only to look up to.

\end{document}